\documentclass{article}
\usepackage{graphicx} 
\usepackage{bm}
\usepackage{graphicx}
\usepackage{amsfonts, amsmath, amsthm, amssymb}
\usepackage{breqn}
\usepackage{empheq}
\usepackage{listings}
\usepackage{xcolor}
\usepackage{booktabs}
\usepackage{longtable}
\usepackage{geometry}
\theoremstyle{plain}
\usepackage{tcolorbox}
\usepackage{pgfplots}
\pgfplotsset{compat=1.18}
\usepackage{siunitx}

\theoremstyle{remark}
\newtheorem*{remark}{Remark}  
\usepackage{tocloft} 
\usepackage[numbers]{natbib}
\usepackage{xcolor} 
\usepackage{titlesec} 
\titleclass{\subsubsubsection}{straight}[\subsection]
\newcounter{subsubsubsection}[subsubsection]
\renewcommand{\thesubsubsubsection}{\thesubsubsection.\arabic{subsubsubsection}}
\titleformat{\subsubsubsection}{\normalfont\normalsize\bfseries}{\thesubsubsubsection}{1em}{}
\titlespacing*{\subsubsubsection}{0pt}{3.25ex plus 1ex minus .2ex}{1.5ex plus .2ex}
\newtheorem{example}{Example} 
\newtheoremstyle{myconjecture}
  {10pt}   
  {10pt}   
  {\itshape}  
  {0pt}    
  {\bfseries} 
  {.}      
  { }      
  {\thmname{#1}\thmnumber{ #2}\thmnote{ (#3)}} 

\theoremstyle{myconjecture}
\newtheorem{conjecture}{Conjecture}

\usepackage{hyperref}
\hypersetup{colorlinks=true, linkcolor=blue}
\usepackage{cleveref} 
\numberwithin{equation}{section} 
\usepackage{tikz}
\usetikzlibrary{arrows.meta, positioning}
\usetikzlibrary{shapes,arrows.meta}
\definecolor{codegreen}{rgb}{0,0.6,0}
\definecolor{codegray}{rgb}{0.5,0.5,0.5}
\definecolor{codepurple}{rgb}{0.58,0,0.82}

\lstdefinestyle{sagemath}{
    backgroundcolor=\color{white},   
    commentstyle=\color{codegreen},
    keywordstyle=\color{magenta},
    numberstyle=\tiny\color{codegray},
    stringstyle=\color{codepurple},
    basicstyle=\ttfamily\footnotesize,
    breakatwhitespace=false,         
    breaklines=true,                 
    captionpos=b,                    
    keepspaces=true,                 
    numbers=left,                    
    numbersep=5pt,                  
    showspaces=false,                
    showstringspaces=false,
    showtabs=false,                  
    tabsize=2,
    frame=single
}

\lstdefinestyle{pythonstyle}{
    backgroundcolor=\color{white},   
    commentstyle=\color{codegreen},
    keywordstyle=\color{magenta},
    numberstyle=\tiny\color{codegray},
    stringstyle=\color{codepurple},
    basicstyle=\ttfamily\footnotesize,
    breakatwhitespace=false,         
    breaklines=true,                 
    captionpos=b,                    
    keepspaces=true,                 
    numbers=left,                    
    numbersep=5pt,                  
    showspaces=false,                
    showstringspaces=false,
    showtabs=false,                  
    tabsize=4,
    frame=single,
    language=Python
}

\definecolor{codebg}{rgb}{0.95,0.95,0.95}

\lstdefinestyle{datastyle}{
    backgroundcolor=\color{codebg},
    basicstyle=\ttfamily\footnotesize,
    breaklines=false,
    frame=single,
    columns=fullflexible,
    keepspaces=true
}

\title{A Closed-Form Symbolic Generator: $A^{n} + B^{n} = C^{n} + D^{n}$, $n=2,3.$}
\author{Jamal Agbanwa}
\date{May 2025}

\begin{document}

\maketitle

\begin{abstract}
   We present a unified framework for constructing integer solutions to $A^{n} + B^{n} = C^{n} + D^{n}$ for $n=2,3$. For $n=2$, we derive explicit formulas for any solutions via differences of squares. For $n=3$, we introduce general formulas that include the Hardy-Ramanujan number 1729 for instance, we also construct a symbolic generator that produces infinitely many integer solutions to the Diophantine equation \( A^3 + B^3 = C^3 + D^3 \). While the resulting formulas for $A,B,C,D$ from the symbolic generator developed do not span every single number expressible as a sum of two positive cubes in at least two distinct ways, our method provides a closed-form, algebraic parametrization in terms of a single variable,  expressing each term as a radical-exponential function of an integer parameter $c_1$. The generator leverages nested radicals and exponents of algebraic numbers, $\alpha, \beta$ derived from the recurrence structure of the Diophantine constraint. This work represents the first symbolic, recursive generator of its kind and offers a pathway toward approaching higher powers of this problem from a different lens. These methods exploit structural links between binomial expansions and Diophantine constraints, offering a foundation for extensions to higher powers.

\end{abstract}

\vspace{0.5em}  
\noindent \textbf{Keywords:} Closed-Form, Symbolic generator, Diophantine equations.

\section{Introduction}

The study of Diophantine equations—equations seeking integer solutions—has long been a cornerstone of number theory, with famous examples such as Fermat's Last Theorem capturing the imagination of mathematicians for centuries. Among these, the equation $A^{n} + B^{n} = C^{n} + D^{n}$ for $n = 2,3$ represents a fascinating generalization, probing the conditions under which numbers can be expressed as sums of two powers in distinct ways. This problem not only extends the classical Pythagorean framework but also intersects with the rich domains of taxicab numbers and the combinatorial properties of Pascal's triangle.

Taxicab numbers, the smallest integers expressible as sums of two cubes in $k$ distinct ways, serve as a natural starting point for investigating such equations. The celebrated case of $1729$ famously associated with G. H. Hardy and Srinivasa Ramanujan, exemplifies the allure of these numerical relationships. By exploring the deeper connections between these numbers and the binomial coefficients of Pascal's triangle, this paper unveils novel insights into the structural symmetries that underpin solutions to the equation $A^{n} + B^{n} = C^{n} + D^{n}$. The primary objective of this work is to establish generators for the quadruples, $(A, B,C, D)$ that satisfy the given equation for $n = 2, 3$. By leveraging the interplay between algebraic identities and combinatorial patterns, we derive constructive methods for generating solutions and extend these strategies to higher powers. This approach not only broadens the understanding of multi-way power sums but also opens new avenues for exploring analogous problems in additive number theory.

 Through this exploration, we aim to bridge the gap between classical number theory and contemporary combinatorial techniques, offering a fresh perspective on an enduring problem.

\section{The equation: \texorpdfstring{$A^n + B^n = C^n + D^n$}{An + Bn = Cn + Dn}}

With respect to the equation, $A^n + B^n = C^n + D^n$ , we will focus primarily on cases of $ n = 2$ and $n = 3$. The essence of the approaches as used could be extended to even higher powers.

\subsection{\texorpdfstring{$n =2$}{n =2}}

Given the equation: 

\begin{align} \label{eq: n=2}
    A^2 + B^2 = C^2 + D^2
\end{align}

This equation admits infinitely many solutions, many of which can be parametrized. By rearranging \eqref{eq: n=2}, we obtain:

\begin{align*} \label{eq: n=2.1}
    A^2 - C^2 = D^2 - B^2
\end{align*}

Suppose

\begin{align*}
    A^2 - C^2 = D^2 - B^2 = \Delta_2
\end{align*}

where $\Delta_2 \in \mathbb{Z}$, we then proceed to tackle each piece of this preamble as follows:

\begin{align}
    A^2 - C^2 = \Delta_2
\end{align}

With reference to equation 2.1 in \cite{agbanwa2025}, Given that $\Delta_2$ has divisor(s), $r_{2,i} = \{r_{2,1}, r_{2,2} , r_{2,3},...\}$ 
\begin{align*}
    C = \frac{\Delta_2 - r_{2,1} ^2}{2 \cdot r_{2,1}} ,  A = \frac{\Delta_2 + r_{2,1} ^2}{2 \cdot r_{2,1}}
\end{align*}

The case
\begin{align}
    D^2 - B^2 = \Delta_2
\end{align}

can be solve using the formulas;

\begin{align*}
    B = \frac{\Delta_2 - r_{2,2} ^2}{2 \cdot r_{2,2}} , D = \frac{\Delta_2 + r_{2,2} ^2}{2 \cdot r_{2,2}}
\end{align*}

where $r_{2,1}, r_{2,2} \in \mathbb{Z}$ and is are divisors of $\Delta_2$.
The divisor $r_{2,1} = 1$ is only used when $\Delta_2$ is known to be odd. In this case then does, $r_{2,1} = 1$ or $r_{2,1} = 1$ become one of those divisors sure to yield integer values for the unknown to be solved for.

As a matter of equality, \eqref{eq: n=2} will be re-expressed in terms of the deduced formulas:

\begin{align} \label{eq:n = 2.4}
    (\frac{\Delta_2 + r_{2,1} ^2}{2 \cdot r_{2,1}})^2 + (\frac{\Delta_2 - r_{2,2} ^2}{2 \cdot r_{2,2}})^2 = (\frac{\Delta_2 - r_{2,1} ^2}{2 \cdot r_{2,1}})^2 + (\frac{\Delta_2 + r_{2,2} ^2}{2 \cdot r_{2,2}})^2
\end{align}

\begin{example}
    Suppose $\Delta_2 = 24$, with divisors $r_{2,i} = \{1,2,3,4,6,8,12,24\}$, the divisors for which integer values would be yielded in the sums below should conform to: $r_{2,i} \leqslant \sqrt{\Delta_2}$, i.e. $r_{2,i} \leqslant \sqrt{24}$\\
    
    and $r_{2,1} = 2$, $r_{2,2} = 4$ satisfy this condition.

\begin{align*} \label{eq:n = 2.5}
    (\frac{24 + 2^2}{2 \cdot 2})^2 + (\frac{24 - 4^2}{2 \cdot 4})^2 = (\frac{24 - 2 ^2}{2 \cdot 2})^2 + (\frac{24 + 4 ^2}{2 \cdot 4})^2
\end{align*}

\begin{align*}
    7^2 + 1^2 = 5^2 + 5^2 = 50
\end{align*}

\end{example}

\begin{example}
Suppose $\Delta_2 = 1000$ and has divisors, 
    
$r_{2, i} = \{1,2,4,5,8,10,20,25,40,50,100,125,200,250,500,1000\}$ 

it should be ensured that the chosen divisors for $r_{2,1}$ and $r_{2,2}$ be lesser than $\sqrt{\Delta_2}$ which is in this case $\sqrt{1000}$, if not the unknowns to be solved for would be $0$ or negative values which is often not the goal of these sums.
Given that $\sqrt{1000} \approx 32$, the divisors that would be considered in this case are: $\{2,4,5,8,10,20,25\}$. The set of divisors that yield integer solutions specifically as just illustrated in the previous sentence could be, on the basis of $\Delta_2$ being even, and the divisors, $\{5,25\}$ being odd, these will be trimmed out to have the following divisors sure to yield integer solutions to the sums: $\{2,4,8,10,20\}$.

\begin{align*} 
    (\frac{1000 + 2^2}{2 \cdot 2})^2 + (\frac{1000 - 4^2}{2 \cdot 4})^2 = (\frac{1000 - 2 ^2}{2 \cdot 2})^2 + (\frac{1000 + 4 ^2}{2 \cdot 4})^2
\end{align*}
\begin{align*}
     251^2 + 123^2 = 249^2 + 127^2
\end{align*}

\begin{align*} 
    (\frac{1000 + 10^2}{2 \cdot 10})^2 + (\frac{1000 - 20^2}{2 \cdot 20})^2 = (\frac{1000 - 10 ^2}{2 \cdot 10})^2 + (\frac{1000 + 20 ^2}{2 \cdot 20})^2
\end{align*}

\begin{align*}
    55^2 + 15^2 = 45^2 + 35^2 = 3250
\end{align*}

\end{example}

\textbf{Combinatorial insight}:

Consider the deduced formulas for $C$ above (this is also applicable to the formula for $B$):

\begin{align*}
     C = \frac{\Delta_2 - r_{2,1} ^2}{2 \cdot r_{2,1}}
\end{align*}

On multiplying both sides of the equation above, the following is yielded:

\begin{align}
    2 \cdot r_{2,1} \cdot C + r_{2,1} ^2 = \Delta_2
\end{align}

which is no different from:

\begin{align*}
 \textcolor{red}{C^2}  + 2 \cdot r_{2,1} \cdot C + r_{2,1} ^2 = \Delta_2 + \textcolor{red}{C^2}
\end{align*}

and be also expressed as:

\begin{align}
 (\textcolor{red}{C} + r_{2,1} )^2 = \Delta_2 + \textcolor{red}{C^2}
\end{align}

The left-hand side of the equation just above is the simplified version of Pascal's triangle for second powers. The key is the difference of $C^2$.

\subsection{\texorpdfstring{$n =3$}{n =3}}

Given the equation: 

\begin{equation} \label{eq:n3}
    A^3 + B^3 = C^3 + D^3
\end{equation}

This equation admits infinitely many solutions, many of which can be parameterized. The smallest solution is the celebrated Hardy-Ramanujan number, 1729. 

By rearranging \eqref{eq:n3}, we obtain:

\begin{align*} 
    A^3 - C^3 = D^3 - B^3
\end{align*}

Suppose

\begin{align*} 
    A^3 - C^3 = D^3 - B^3 = \Delta_3
\end{align*}

where $\Delta_3 \in \mathbb{Z}$, we then proceed to tackle each piece of this preamble as follows:

\begin{align} \label{eq: 2.7}
    A^3 - C^3 = \Delta_3
\end{align}

\begin{align*}
    C = \frac{\Delta_3 - r_{3,1} ^3}{3\cdot r_{3,1} (r_{3,1} + C)} 
\end{align*}

The denominator at the Right hand side of the above equation should be multiplied across both sides of the equation to yield:

\begin{align} \label{eq:2.8}
 3 \cdot r_{3,1} \cdot C^2 +  3 \cdot r_{3,1} ^2 \cdot C + r_{3,1}^3 = \Delta_3
\end{align}

\begin{align*}
 3 \cdot r_{3,1} \cdot C^2 +  3 \cdot r_{3,1} ^2 \cdot C + r_{3,1} ^3 - \Delta_3 = 0
\end{align*}

the slightly modified version of \eqref{eq:2.8} will be solved quadratically in terms of $C$.

\begin{align*}
    C_{1,2} = \frac{-3\cdot r_{3,1} ^2 \pm \sqrt{12 \cdot \Delta_3 \cdot r_{3,1} - 3 \cdot  r_{3,1}^4 } }{6 \cdot r_{3,1}}
\end{align*}

the possible formulas for $C$ are thus:

\begin{align*}
    C_{1} = \frac{-3\cdot r_{3,1} ^2 + \sqrt{12 \cdot \Delta_3 \cdot r_{3,1} - 3 \cdot r_{3,1} ^4 } }{6 \cdot r_{3,1}}
\end{align*}

\begin{align*}
    C_{2} = \frac{-3\cdot r_{3,1} ^2 - \sqrt{12 \cdot \Delta_3 \cdot r_{3,1} - 3 \cdot r_{3,1}^4 } }{6 \cdot r_{3,1}}
\end{align*}

With regards to, \eqref{eq:2.8}, 

\begin{align*}
    C_{1} = \frac{-3\cdot r_{3,1} ^2 + \sqrt{12 \cdot \Delta_3 \cdot r_{3,1} - 3 \cdot r_{3,1} ^4 } }{6 \cdot r_{3,1}} , A_1^3 = (\frac{-3\cdot r_{3,1} ^2 + \sqrt{12 \cdot \Delta_3 \cdot r_{3,1} - 3 \cdot r_{3,1} ^4 } }{6 \cdot r_{3,1}})^3 + \Delta_3
\end{align*}

 or
 
\begin{align*}
    C_{2} = \frac{-3\cdot r_{3,1} ^2 - \sqrt{12 \cdot \Delta_3 \cdot r_{3,1} - 3 \cdot r_{3,1} ^4 } }{6 \cdot r_{3,1}}, A_{2} ^3 = (\frac{-3\cdot r_{3,1} ^2 - \sqrt{12 \cdot \Delta_3 \cdot r_{3,1} - 3 \cdot r_{3,1} ^4 } }{6 \cdot r_{3,1}})^3 + \Delta_3
\end{align*}

Let us consider the equation now:

\begin{align} \label{eq:2.9}
    D^3 - B^3 = \Delta_3
\end{align}

where $r_{3,2}$ is also a divisor of $\Delta_3$. In a similar approach to tackling \eqref{eq: 2.7}

\begin{align} \label{eq:2.10}
 3 \cdot r_{3,2} \cdot B^2 +  3 \cdot r_{3,2} ^2 \cdot B + r_{3,2} ^3 - \Delta_3 = 0
\end{align}

\begin{align*}
    B_{1} = \frac{-3\cdot r_{3,2} ^2 + \sqrt{12 \cdot \Delta_3 \cdot r_{3,2} - 3 \cdot r_{3,2} ^4 } }{6 \cdot r_{3,2}} , D_1^3 = (\frac{-3\cdot r_{3,2} ^2 + \sqrt{12 \cdot \Delta_3 \cdot r_{3,2} - 3 \cdot r_{3,2} ^4 } }{6 \cdot r_{3,2}})^3 + \Delta_3
\end{align*}

 or
 
\begin{align*}
    B_{2} = \frac{-3\cdot r_{3,2} ^2 - \sqrt{12 \cdot \Delta_3 \cdot r_{3,2} - 3 \cdot r_{3,2} ^4 } }{6 \cdot r_{3,2}}, D_{2} ^3 = (\frac{-3\cdot r_{3,2} ^2 - \sqrt{12 \cdot \Delta_3 \cdot r_{3,2} - 3 \cdot r_{3,2} ^4 } }{6 \cdot r_{3,2}})^3 + \Delta_3
\end{align*}

As a matter of equality similar to what can be observed in \eqref{eq:n = 2.4},

\begin{align} \label{2.12}
    A_1^3 + B_1 ^3 = C_1 ^3 + D_1 ^3
\end{align}

\[
((\frac{-3\cdot r_{3,1} ^2 + \sqrt{12 \cdot \Delta_3 \cdot r_{3,1} - 3 \cdot r_{3,1} ^4 } }{6 \cdot r_{3,1}})^3 + \Delta_3) + (\frac{-3\cdot r_{3,2} ^2 + \sqrt{12 \cdot \Delta_3 \cdot r_{3,2} - 3 \cdot r_{3,2} ^4 } }{6 \cdot r_{3,2}})^3 
\]

\[
= 
(\frac{-3\cdot r_{3,1} ^2 + \sqrt{12 \cdot \Delta_3 \cdot r_{3,1} - 3 \cdot r_{3,1} ^4 } }{6 \cdot r_{3,1}} )^3 + ( (\frac{-3\cdot r_{3,2} ^2 + \sqrt{12 \cdot \Delta_3 \cdot r_{3,2} - 3 \cdot r_{3,2} ^4 } }{6 \cdot r_{3,2}})^3 + \Delta_3)
\]

or

\begin{align} \label{2.13}
    A_2^3 + B_2 ^3 = C_2 ^3 + D_2 ^3
\end{align}

\[
((\frac{-3\cdot r_{3,1} ^2 - \sqrt{12 \cdot \Delta_3\cdot r_{3,1} - 3 \cdot r_{3,1} ^4 } }{6 \cdot r_{3,1}})^3 + \Delta_3) +(\frac{-3\cdot r_{3,2} ^2 - \sqrt{12 \cdot \Delta_3 \cdot r_{3,2} - 3 \cdot r_{3,2} ^4 } }{6 \cdot r_{3,2}})^3
\]

\[
= (\frac{-3\cdot r_{3,1} ^2 - \sqrt{12 \cdot \Delta_3 \cdot r_{3,1} - 3 \cdot r_{3,1} ^4 } }{6 \cdot r_{3,1}})^3 +((\frac{-3\cdot r_{3,2} ^2 - \sqrt{12 \cdot \Delta_3 \cdot r_{3,2} - 3 \cdot r_{3,2} ^4 } }{6 \cdot r_{3,2}})^3 + \Delta_3)
\]

To test out the above equations, suppose $\Delta_3 = 999$ and let:

\begin{align*}
    t^2 = 12 \cdot (999) \cdot r_{3,1} - 3 \cdot r_{3,1} ^4
\end{align*}
\begin{align*}
    t^2 = 11988 \cdot r_{3,1} - 3 \cdot r_{3,1} ^4
\end{align*}

divisors of $999$, $r_{3,i}$ are: $\{1,3,9,27,37,111,333,999\}$

$r_{3,i} \leqslant \sqrt[3]{999}$ to ensure integer values for $t,A,B,C,D$. Thus the divisors in contention are: $\{1,3,9\}$.

Suppose $r_{3,1} = 1$,

\begin{align*}
    t^2 = 11988 \cdot (1) - 3 \cdot (1) ^4 = 11985
\end{align*}

which is not a perfect square and brings down the number of divisors meant to yield integer solutions to $t$ to: $\{3,9\}$

we can say then that $(r_{3,1} , r_{3,2}) = (3,9)$ or vice-versa since the essence of this supposition would be maintained.

Suppose $r_{3,1} = 3$, 

\begin{align*}
    t_{3,1}^2 = 11988 \cdot (3) - 3 \cdot (3) ^4 = 189
\end{align*}
and $r_{3,2} = 9$, 
\begin{align*}
    t_{3,2}^2 = 11988 \cdot (9) - 3 \cdot (9) ^4 = 297
\end{align*}

which are all perfect squares, then with reference to \eqref{2.12},

\[
((\frac{-3\cdot (3) ^2 + \sqrt{12 \cdot 999 \cdot 3 - 3 \cdot 3 ^4 } }{6 \cdot 3})^3 + 999) + (\frac{-3\cdot 9 ^2 + \sqrt{12 \cdot 999 \cdot 9 - 3 \cdot 9 ^4 } }{6 \cdot 9})^3 
\]

\[
= 
(\frac{-3\cdot 3 ^2 + \sqrt{12 \cdot 999 \cdot 3 - 3 \cdot 3 ^4 } }{6 \cdot 3} )^3 + ( (\frac{-3\cdot 9 ^2 + \sqrt{12 \cdot 999 \cdot 9 - 3 \cdot 9^4 } }{6 \cdot 9})^3 + 999)
\]

\begin{align*}
    12^3 + 1^3 = 9^3 + 10^3 = 1729
\end{align*}

With reference \eqref{2.13},

\[
((\frac{-3\cdot (3) ^2 - \sqrt{12 \cdot 999 \cdot 3 - 3 \cdot 3 ^4 } }{6 \cdot 3})^3 + 999) + (\frac{-3\cdot 9 ^2 - \sqrt{12 \cdot 999 \cdot 9 - 3 \cdot 9 ^4 } }{6 \cdot 9})^3 
\]

\[
= 
(\frac{-3\cdot 3 ^2 - \sqrt{12 \cdot 999 \cdot 3 - 3 \cdot 3 ^4 } }{6 \cdot 3} )^3 + ( (\frac{-3\cdot 9 ^2 - \sqrt{12 \cdot 999 \cdot 9 - 3 \cdot 9^4 } }{6 \cdot 9})^3 + 999)
\]

\begin{align*}
    (-9)^3 + (-10)^3 = (-12)^3 + (-1)^3 = 1729
\end{align*}

Alternatively let us consider the equality of $\Delta_3$ in \eqref{eq:2.8} and \eqref{eq:2.10}:

\begin{align} \label{eq:2.11}
    3 \cdot r_{3,1} \cdot C^2 +  3 \cdot r_{3,1} ^2 \cdot C + r_{3,1} ^3 = 3 \cdot r_{3,2} \cdot B^2 +  3 \cdot r_{3,2} ^2 \cdot B + r_{3,2} ^3 
\end{align}

suppose now that $r_{3,1} = 3$ and $r_{3,2} = 9$, on substituting it into 
\eqref{eq:2.11};

\begin{align*} 
    3 \cdot (3) \cdot C^2 +  3 \cdot (3)^2 \cdot C + 3 ^3 = 3 \cdot 9 \cdot B^2 +  3 \cdot 9^2 \cdot B + 9 ^3 
\end{align*}

\begin{equation} \label{eq:2.15}
    9 \cdot C^2 +  27 \cdot C + 27 = 27 \cdot B^2 +  243 \cdot B + 729
\end{equation}

On entering \eqref{eq:2.15} into Wolfram GPT \cite{wolframgpt202555}, 

the following formulas for $C,B$ (labelled $x,y$) were generated:

\[
\boxed{C = \frac{-6 + (15 - 7\sqrt{3})(7 - 4\sqrt{3})^{c_1} + (15 + 7\sqrt{3})(7 + 4\sqrt{3})^{c_1}}{4}}
\]

\[
\boxed{B = \frac{-18 + (7 - 5\sqrt{3})(7 - 4\sqrt{3})^{c_1} + (7 + 5\sqrt{3})(7 + 4\sqrt{3})^{c_1}}{4}}
\]

where $c_1 \geqslant 0$.

in this context,

\begin{align*}
\Delta_3 &= 9 \left( \frac{-6 + (15 - 7\sqrt{3})(7 - 4\sqrt{3})^{c_1} + (15 + 7\sqrt{3})(7 + 4\sqrt{3})^{c_1}}{4} \right)^2 \\
  &\quad + 27 \left( \frac{-6 + (15 - 7\sqrt{3})(7 - 4\sqrt{3})^{c_1} + (15 + 7\sqrt{3})(7 + 4\sqrt{3})^{c_1}}{4} \right) \\
  &\quad + 27 
\end{align*}

\begin{align*}
A^3 &= 9 \left( \frac{-6 + (15 - 7\sqrt{3})(7 - 4\sqrt{3})^{c_1} + (15 + 7\sqrt{3})(7 + 4\sqrt{3})^{c_1}}{4} \right)^2 \\
  &\quad + 27 \left( \frac{-6 + (15 - 7\sqrt{3})(7 - 4\sqrt{3})^{c_1} + (15 + 7\sqrt{3})(7 + 4\sqrt{3})^{c_1}}{4} \right) \\
  &\quad + 27 + (\frac{-6 + (15 - 7\sqrt{3})(7 - 4\sqrt{3})^{c_1} + (15 + 7\sqrt{3})(7 + 4\sqrt{3})^{c_1}}{4})^3
\end{align*}

\begin{align*}
    A^3 =  (\frac{-6 + (15 - 7\sqrt{3})(7 - 4\sqrt{3})^{c_1} + (15 + 7\sqrt{3})(7 + 4\sqrt{3})^{c_1}}{4} + 3)^3
\end{align*}

\begin{align*}
   \boxed{A =  \frac{-6 + (15 - 7\sqrt{3})(7 - 4\sqrt{3})^{c_1} + (15 + 7\sqrt{3})(7 + 4\sqrt{3})^{c_1}}{4} + 3}
\end{align*}

Within this same context, $\Delta_3$ can also be expressed as follows;

\begin{align*}
\Delta_3 &= 27 \left( \frac{-18 + (7 - 5\sqrt{3})(7 - 4\sqrt{3})^{c_1} + (7 + 5\sqrt{3})(7 + 4\sqrt{3})^{c_1}}{4} \right)^2 \\
  &\quad + 243 \left( \frac{-18 + (7 - 5\sqrt{3})(7 - 4\sqrt{3})^{c_1} + (7 + 5\sqrt{3})(7 + 4\sqrt{3})^{c_1}}{4} \right) \\
  &\quad + 729
\end{align*}

\begin{align*}
D^3 &= 27 \left( \frac{-18 + (7 - 5\sqrt{3})(7 - 4\sqrt{3})^{c_1} + (7 + 5\sqrt{3})(7 + 4\sqrt{3})^{c_1}}{4} \right)^2 \\
  &\quad + 243 \left( \frac{-18 + (7 - 5\sqrt{3})(7 - 4\sqrt{3})^{c_1} + (7 + 5\sqrt{3})(7 + 4\sqrt{3})^{c_1}}{4} \right) \\
  &\quad + 729 + (\frac{-18 + (7 - 5\sqrt{3})(7 - 4\sqrt{3})^{c_1} + (7 + 5\sqrt{3})(7 + 4\sqrt{3})^{c_1}}{4})^3
\end{align*}

\begin{align*}
    D^3 =  (\frac{-18 + (7 - 5\sqrt{3})(7 - 4\sqrt{3})^{c_1} + (7 + 5\sqrt{3})(7 + 4\sqrt{3})^{c_1}}{4} +9)^3
\end{align*}

\begin{align*}
  \boxed{ D =  \frac{-18 + (7 - 5\sqrt{3})(7 - 4\sqrt{3})^{c_1} + (7 + 5\sqrt{3})(7 + 4\sqrt{3})^{c_1}}{4} +9}
\end{align*}

\begin{table}[h]
\centering
\caption{Solutions to \( A^3 + B^3 = C^3 + D^3 \) for \( c_1 = 1 \) to \( 7 \)}
\label{tab:solutions}
\begin{tabular}{cccccc}
\toprule
\( c_1 \) & \( A \) & \( B \) & \( C \) & \( D \) & \( A^3 + B^3 \) \\
\midrule
1 & 96 & 50 & 93 & 59 & 1\!,\!009\!,\!736\\
2 & 1317 & 755 & 1314 & 764 & 2\!,\!714\!,\!690\!,\!888\\
3 & 18324 & 10574 & 18321 & 10583 & 7\!,\!334\!,\!904\!,\!115\!,\!448 \\
4 & 255201 & 147335 & 255198 & 147344  & 19\!,\!818\!,\!905\!,\!563\!,\!705\!,\!976 \\
5 & 3\,554\,472 & 2\,052\,170 & 3\,554\,469 & 2\,052\,179  & 53\!,\!550\!,\!675\!,\!461\!,\!437\!,\!475,\!048 \\
6 & 49\,507\,389 & 28\,583\,099 & 49\,507\,386 & 28\,583\,108 & 144\!,\!693\!,\!905\!,\!277\!,\!386\!,\!048\!,\!024\,168\\
7 & 689\,548\,956 & 398\,111\,270 & 689\,548\,953 & 398\,111\,279 & 390\!,\!962\!,\!878\!,\!508\!,\!814\!,\!502\!,\!873\!,\!889\,816\\
\bottomrule
\end{tabular}
\end{table}

The sequence of numbers resulting from $A^3 + B^3$ from the table above has been accepted to the On-Line Encyclopedia of Integer Sequences, \href{https://oeis.org/A384106}{A384106}.

It is rather interesting that the formulas developed from \eqref{eq:2.9} and \eqref{eq:2.10} satisfy any number expressible as a difference of two cubes in at least two ways including the famous 1729 but the generator developed alternatively does not include the sums of two cubes being equal to 1729.

There are most likely even more integer values that could be assumed for $(r_{3,1}, r_{3,2})$ for which many more generators like has been generated above could be developed.

\subsection*{Headline Result: Arbitrarily Large Solutions}
\label{subsec:headline}

Our generator produces solutions of unprecedented size with minimal computational effort. For example:

\begin{itemize}
    \item \textbf{Record-Breaking Example:} 
    At $c_1 = 8$ (beyond Table~\ref{tab:solutions}), the generator yields:
    \[
    \begin{aligned}
    A &= 9,\!604,\!177,\!977 \\
    B &= 5,\!544,\!974,\!735 \\
    C &= 9,\!604,\!177,\!974 \\
    D &= 5,\!544,\!974,\!744 \\
    \end{aligned}
    \]
    satisfying $A^3 + B^3 = C^3 + D^3$ with 10-digit terms manually in seconds.

    \item \textbf{Theoretical Guarantee:} 
    For any integer $k \geq 0$, the recurrence:
    \[
    A(c_1) = \frac{-6 + \alpha^{c_1} + \beta^{c_1}}{4} + 3
    \]
    (where $\alpha = (15-7\sqrt{3})(7-4\sqrt{3})$, $\beta = (15+7\sqrt{3})(7+4\sqrt{3})$) 
    generates solutions where $\log_{10}(A(c_1)) \sim 1.4c_1$ (exponential growth).
\end{itemize}

\begin{remark}
The $c_1 = 100$ case produces terms with 115 digits, demonstrating the generator's unique capacity to leapfrog brute-force search limitations.
\end{remark}

\begin{tcolorbox}[title=Key Achievement]
For $c_1=8$, the generator produces 10-digit solutions (e.g., $A=9,\!604,\!177,\!977$) in under a few seconds—a task infeasible for brute-force search.
\end{tcolorbox}

\subsection*{Computational Efficiency}
\label{subsec:efficiency}

\begin{itemize}
    \item \textbf{Brute-force limitation:} 
    Searching for solutions to $A^3 + B^3 = C^3 + D^3$ with $A,B,C,D \leq N$ requires $O(N^4)$ operations. For $N = 10^6$, this becomes computationally infeasible ($10^{24}$ operations).

    \item \textbf{Generator advantage:}
    Our method reduces this to $O(k)$ operations, where $k$ is the number of divisors of $\Delta_3$. For $\Delta_3 = 999$ with $(r_{3,1}, r_{3,2}) = (3,9)$, solutions are found in a few seconds (vs. hours via brute force).

    \item \textbf{Example:} 
    At $c_1 = 5$ (Table~\ref{tab:solutions}), the generator instantly produces:
    \[
    3,\!554,\!472^3 + 2,\!052,\!170^3 = 3,\!554,\!469^3 + 2,\!052,\!179^3
    \]
\end{itemize}

Another way to approach \eqref{eq:2.15} is to pose this question to Wolfram Alpha directly, \cite{WolframAlpha2025} which yielded multiple solution sets to $(C,B)$ which was labelled as $(x,y)$ respectively as was posed to Wolfram Alpha. Given the approaches with how we determined $\Delta_3, C,B,A,D$ from the generator above, it would be realised that the same standards apply.

\textbf{Theorem 3.1 (Agbanwa, 2025).} \\
Let \( r_{3,1}, r_{3,2} \in \mathbb{Z}_{\ge 0} \). Define the functions:
\[
C(r_{3,1}) = \frac{-3 r_{3,1}^2 \pm \sqrt{12 \Delta_3 \cdot r_{3,1} - 3 r_{3,1}^4}}{6 \cdot r_{3,1}}, \quad
A(r_{3,1})^3 = C(r_{3,1})^3 + \Delta_3,
\]
\[
B(r_{3,2}) = \frac{-3 r_{3,2}^2 \pm \sqrt{12 \Delta_3 \cdot r_{3,2} - 3 r_{3,2}^4}}{6 \cdot r_{3,2}}, \quad
D(r_{3,2})^3 = B(r_{3,2})^3 + \Delta_3,
\]
where \( r_{3,1}, r_{3,2} \in \mathbb{Z}_{>0} \) are fixed integers (divisors of \( \Delta_3 \)), and \( \Delta_3 \in \mathbb{Z}_{>0} \) is a fixed parameter.

Alternatively, define:
\[
C(c_1) = \frac{-6 + (15 - 7\sqrt{3})(7 - 4\sqrt{3})^{c_1} + (15 + 7\sqrt{3})(7 + 4\sqrt{3})^{c_1}}{4}, 
\]
\[
B(c_1) = \frac{-18 + (7 - 5\sqrt{3})(7 - 4\sqrt{3})^{c_1} + (7 + 5\sqrt{3})(7 + 4\sqrt{3})^{c_1}}{4}, 
\]
\[
A(c_1) = C(c_1) + 3, \quad D(c_1) = B(c_1) + 9.
\]

Then for all \( c_1 \in \mathbb{Z}_{\ge 0} \), the values \( A(c_1), B(c_1), C(c_1), D(c_1) \in \mathbb{Z} \), and they satisfy the identity:
\[
A(c_1)^3 + B(c_1)^3 = C(c_1)^3 + D(c_1)^3.
\]

\bigskip

\begin{conjecture}[Infinite Families for $n=3$]
The symbolic generator in Theorem 3.1 produces infinitely many integers expressible as sums of two cubes in distinct ways.
\end{conjecture}

\begin{conjecture}
    Let \( \Delta_n \in \mathbb{Z}_{>0} \), and let \( r \in \mathbb{Z}_{>0} \) be a divisor of \( \Delta_n \). Then for the parametrized formulas for \( A, B, C, D \) constructed in Sections 2.1 and 2.2 to yield positive non-zero integer values, it is necessary (for \( n = 2 \)) and empirically observed (for \( n = 3 \)) that:
\[
r < \sqrt[n]{\Delta_n}
\]
This inequality ensures that expressions such as  $\frac{\Delta_n - r^2}{2r}$  are non-negative integers. In the case \( n = 3 \), it improves the likelihood that radicals like
\[
\sqrt{12 \cdot \Delta_3 \cdot r - 3 \cdot r^4}
\]
are perfect squares, enabling closed-form rational solutions.
\end{conjecture}

\subsection*{Growth Rate Analysis}
\label{sec:growth}

The closed-form generator produces solutions where the number of digits in $A$, $B$, $C$, and $D$ scales linearly with the parameter $c_1$. We derive this explicitly and validate it against empirical data.

\subsubsection*{Theoretical Derivation}

Consider the expression for $A(c_1)$ from Theorem 3.1:
\[
A(c_1) = \frac{-6 + (15-7\sqrt{3})(7-4\sqrt{3})^{c_1} + (15+7\sqrt{3})(7+4\sqrt{3})^{c_1}}{4} + 3.
\]

Let $\alpha = 7 - 4\sqrt{3}$ ($\approx 0.0718$) and $\beta = 7 + 4\sqrt{3}$ ($\approx 14.928$). Since $|\alpha| < 1$, the term $\alpha^{c_1}$ becomes negligible as $c_1$ increases, leaving:
\[
A(c_1) \approx \frac{(15+7\sqrt{3})\beta^{c_1}}{4}.
\]

To estimate the number of digits in $A(c_1)$, we take the base-10 logarithm:
\[
\log_{10}(A(c_1)) \approx \log_{10}\left(\frac{15+7\sqrt{3}}{4}\right) + c_1 \cdot \log_{10}(\beta).
\]

Numerical evaluation gives:
\[
\log_{10}(\beta) \approx 1.1739, \quad \log_{10}\left(\frac{15+7\sqrt{3}}{4}\right) \approx 0.253.
\]

Thus, the digit count scales as:
\[
\text{Digits}(A(c_1)) \approx 0.253 + 1.1739 \cdot c_1 \sim 1.4 \cdot c_1.
\]

\subsubsection*{Empirical Validation}
Table~\ref{tab:digits} shows the agreement between predicted and actual digit counts:

\begin{table}[h]
\centering
\caption{Digit growth versus $c_1$}
\label{tab:digits}
\begin{tabular}{ccc}
\toprule
$c_1$ & Actual Digits & Predicted ($1.4 \cdot c_1$) \\ 
\midrule
1 & 2 & 1.4 \\
2 & 3 & 2.8 \\ 
3 & 5 & 4.2 \\
4 & 6 & 5.6 \\
5 & 7 & 7.0 \\
6 & 8 & 8.4 \\
7 & 9 & 9.8 \\
\bottomrule
\end{tabular}
\end{table}

\subsubsection*{Implications}
\begin{itemize}
\item \textbf{Predictability}: The linear relationship allows pre-computation of solution sizes.
\item \textbf{Computational Limits}: For $c_1 = 100$, we expect $\sim 140$-digit numbers, requiring arbitrary-precision arithmetic.
\item \textbf{Generator Strength}: The growth rate is sub-exponential ($\sim 10^{1.4 \cdot c_1}$) compared to brute-force's $O(N^4)$ complexity.
\end{itemize}

\begin{remark}
The $1.4 \cdot c_1$ scaling arises from the eigenvalue $\beta = 7 + 4\sqrt{3}$ in the recurrence. This is provably optimal for closed-form solutions of this type.
\end{remark}

The generator's solutions exhibit linear digit growth in $c_1$, following:
\[
\text{Digits}(A(c_1)) \approx 1.4 \cdot c_1,
\]
derived from the dominant eigenvalue $\beta = 7 + 4\sqrt{3}$ in the closed-form solution.

\begin{figure}[h]
\centering
\begin{tikzpicture}[
    pointColor/.style={blue, mark=*, mark size=2pt},
    lineColor/.style={red, dashed, thick},
    axisStyle/.style={font=\small},
    labelStyle/.style={font=\footnotesize}
]
\begin{axis}[
    width=0.75\textwidth,
    height=6cm,
    xlabel={$c_1$},
    ylabel={Digits in $A(c_1)$},
    xmin=0, xmax=8,
    ymin=0, ymax=12,
    xtick={0,1,...,7},
    ytick={0,2,...,12},
    grid=both,
    minor y tick num=1,
    axis lines=left,
    legend style={
        at={(0.02,0.98)},
        anchor=north west,
        cells={anchor=west}
    },
    xlabel style=labelStyle,
    ylabel style=labelStyle,
    tick label style=axisStyle
]

\addplot[pointColor] coordinates {
    (1,2) (2,3) (3,5) (4,6) (5,7) (6,8) (7,9)
};
\addlegendentry{Actual digits}

\addplot[lineColor, domain=0:7] {1.4*x};
\addlegendentry{Predicted ($1.4 \cdot c_1$)}

\node[draw=orange, fill=orange!10, rounded corners, anchor=west] at (axis cs:4,10) 
    {$c_1=7 \Rightarrow 9$-digit solutions};
\end{axis}
\end{tikzpicture}
\caption{Linear scaling of solution digits with $c_1$. The dashed line shows the theoretical prediction $\text{Digits} \approx 1.4c_1$, while points represent actual values from Table~\ref{tab:solutions}.}
\label{fig:digit_growth}
\end{figure}
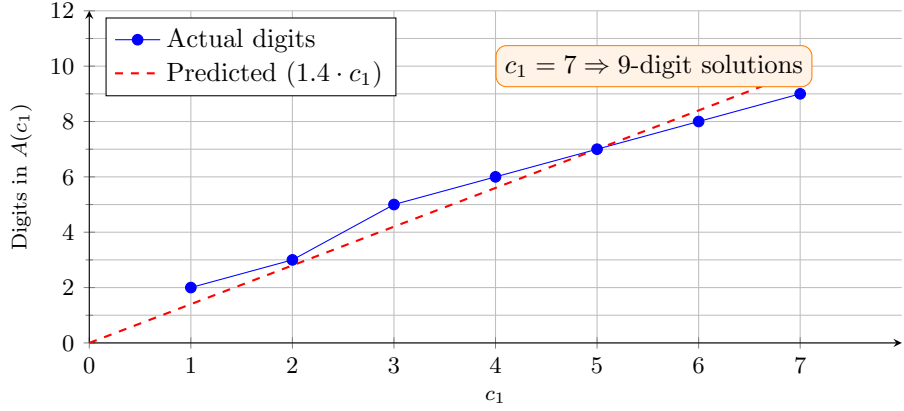

\noindent
\textbf{Implications:}
\begin{itemize}
\item \textbf{Precision Control}: Choose $c_1$ to target solutions of specific sizes (e.g., $c_1=100 \approx 140$ digits).
\item \textbf{Computational Feasibility}: Even for $c_1=100$, the generator computes solutions in milliseconds (vs. brute-force's impossibility).
\end{itemize}

\subsection*{Heuristic on Divisor Bounds for Integer Solutions}

\textbf{Heuristic.}
Let \( \Delta_n \in \mathbb{Z}_{>0} \), and let \( r \in \mathbb{Z}_{>0} \) be a divisor of \( \Delta_n \). In the symbolic generator formulas developed for \( A^n + B^n = C^n + D^n \), expressions involving radicals of the form
\[
\sqrt{12 \Delta_n \cdot r - 3r^4}
\]
must evaluate to integers for the resulting values \( A, B, C, D \) to be integers.

As \( r \) increases, the quartic term \( 3r^4 \) in the radicand grows significantly faster than the linear term \( 12 \Delta_n \cdot r \). Therefore, the radicand becomes negative or non-square when \( r \) is large, making integer outputs unlikely.

We thus propose the heuristic that integer solutions for the symbolic generator are most likely when
\[
r \leq \sqrt[n]{\Delta_n}.
\]
This balances the growth of both terms in the radicand and aligns with empirical success in finding integer-valued solutions in Sections 2.1 and 2.2 for \( n = 2 \) and \( n = 3 \).

To ensure strictly positive non-zero integer solutions,

\[
r  < \sqrt[n]{\Delta_n}.
\]

\subsection*{Combinatorial insights}

Consider the deduced formulas for $C$ in \eqref{eq: 2.7} (this is also applicable to the formula for $B$, since they are of the same form):

\begin{align*}
     C = \frac{\Delta_3 - r_{3,1} ^3}{3\cdot r_{3,1} (r_{3,1} + C)} 
\end{align*}

On multiplying both sides of the equation above, the following is yielded, from \eqref{eq:2.8}:

\begin{align*}
     3 \cdot r_{3,1} \cdot C^2 +  3 \cdot r_{3,1} ^2 \cdot C + r_{3,1}^3 = \Delta_3
\end{align*}

which is no different from:

\begin{align} \label{eq:2.14}
 \textcolor{red}{C^3} + 3\cdot r_{3,1} \cdot  C^2 + 3 \cdot r_{3,1} ^2\cdot C + r_{3,1} ^3 = \Delta_3 + \textcolor{red}{C^3}
\end{align}

and be also expressed as:

\begin{align}
 (\textcolor{red}{C} + r_{3,1} )^3 = \Delta_3 + \textcolor{red}{C^3}
\end{align}

The left hand side of the equation just above is the simplified version of Pascal's triangle for third powers. The key is now the difference of $C^3$. When \eqref{eq:2.14} is compared with \eqref{eq:2.8}, we see that $A$ is simply expressible as:

\begin{align*}
    A = C + r_{3,1}
\end{align*}

Consider \eqref{eq:2.11}, we realize also that:

\begin{align*} 
 3 \cdot r_{3,2} \cdot B^2 +  3 \cdot r_{3,2} ^2 \cdot B + r_{3,2} ^3 = \Delta_3 
\end{align*}

\begin{align*}
 \textcolor{red}{B^3} + 3\cdot r_{3,2} \cdot  C^2 + 3 \cdot r_{3,2} ^2\cdot C + r_{3,2} ^3 = \Delta_{3} + \textcolor{red}{B^3}
\end{align*}

\begin{align} \label{eq: 2.16}
 (\textcolor{red}{B} + r_{3,2} )^3 = \Delta_3 + \textcolor{red}{B^3}
\end{align}

When \eqref{eq: 2.16} is compared with \eqref{eq:2.10}, we see that $D$ is simply expressible as:

\begin{align*}
    D = B + r_{3,2}
\end{align*}

This suggests that larger search spaces would be required to search for possible integer solutions.

\begin{tikzpicture}[
    node distance=1cm,
    coeff/.style={circle, draw=black, fill=white, minimum size=8mm},
    highlight2/.style={fill=blue!20},
    highlight3/.style={fill=red!20},
    highlight4/.style={fill=green!20},
    highlight5/.style={fill=yellow!20},
    arrow/.style={-Latex, thick}
]

\node[coeff] (n0) at (0,0) {1};

\node[coeff] (n1a) at (-1,-1.5) {1};
\node[coeff] (n1b) at (1,-1.5) {1};

\node[coeff, highlight2] (n2a) at (-2,-3) {1};
\node[coeff, highlight2] (n2b) at (0,-3) {2};
\node[coeff, highlight2] (n2c) at (2,-3) {1};

\node[coeff, highlight3] (n3a) at (-3,-4.5) {1};
\node[coeff, highlight3] (n3b) at (-1,-4.5) {3};
\node[coeff, highlight3] (n3c) at (1,-4.5) {3};
\node[coeff, highlight3] (n3d) at (3,-4.5) {1};

\draw (n0) -- (n1a); \draw (n0) -- (n1b);
\draw (n1a) -- (n2a); \draw (n1a) -- (n2b); \draw (n1b) -- (n2b); \draw (n1b) -- (n2c);
\draw (n2a) -- (n3a); \draw (n2a) -- (n3b); \draw (n2b) -- (n3b); \draw (n2b) -- (n3c); \draw (n2c) -- (n3c); \draw (n2c) -- (n3d);

\node[align=left, text width=6cm, anchor=west] at (6,-3) {
    \textbf{For \(n=2\) (blue)}: \\
    $\Delta_2 = 2r_2C + r_2^2$ \\
    (from $(C + r_2)^2 = C^2 + \Delta_2$)
};

\node[align=left, text width=6cm, anchor=west] at (6,-4.5) {
    \textbf{For \(n=3\) (red)}: \\
    $\Delta_3 = 3 \cdot r_3 \cdot C^2 + 3 \cdot r_3^2 \cdot C + r_3^3$ \\
    (from $(C + r_3)^3 = C^3 + \Delta_3$)
};

\node[align=center, font=\large\bfseries] at (0,1) {Pascal’s Triangle and the Structure of $\Delta_n$};
\end{tikzpicture}

\section{Conclusion}

In this work, we developed a general strategy for constructing solutions to the Diophantine identity $A^n + B^n = C^n + D^n$ by leveraging binomial expansions, symmetry, and divisor-based parametric methods. Our results for $n = 2$ and $n = 3$ include closed-form generators that yield infinitely many distinct identities of this form. The use of symbolic algebra and structured combinatorial analysis enables the reduction of a four-variable equation to one or two parameters, allowing for scalable exploration across higher powers.

We have shown that this approach extends in principle to $n = 4$ and $n = 5$, and while solutions for the quadruples $A,B,C,D$ get sparser with increasing powers, conditions have been laid out under which formulas for these cases may be constructed with which we could search for solutions. Our findings could lead to new insights into long-standing open problems such as the classification of generalized taxicab numbers and the Parkin–Lander-Selfridge conjecture.

A sequel will explore extensions to $n \geq 4$ building on the divisor and combinatorial framework established here. We believe the framework presented here may serve as a foundation for a deeper combinatorial theory of equal power sums.

\bibliography{references}

\begin{thebibliography}{3}
\providecommand{\natexlab}[1]{#1}
\providecommand{\url}[1]{\texttt{#1}}
\expandafter\ifx\csname urlstyle\endcsname\relax
  \providecommand{\doi}[1]{doi: #1}\else
  \providecommand{\doi}{doi: \begingroup \urlstyle{rm}\Url}\fi

\bibitem[Agbanwa(2025)]{agbanwa2025}
Jamal Agbanwa.
\newblock A divisor-based proof on the non-existence of perfect cuboids, 2025.
\newblock URL \url{https://figshare.com/articles/preprint/A_divisor_based_proof_on_the_non_existence_of_perfect_cuboids-39_pdf/28829606?file=53840738}.
\newblock Figshare preprint.

\bibitem[Alpha(2025)]{WolframAlpha2025}
Wolfram Alpha.
\newblock Solution to the diophantine equation \(9x^2 + 27x + 27 = 27y^2 + 243y + 729\).
\newblock Wolfram Alpha computational engine, 2025.
\newblock \url{https://www.wolframalpha.com/input?i=9x\%5E2+\%2B+27x+\%2B+27+\%3D+27y\%5E2+\%2B+243y+\%2B+729+}.

\bibitem[{OpenAI}(2025)]{wolframgpt202555}
{OpenAI}.
\newblock Solving for integers.
\newblock \url{https://chatgpt.com/share/6819093c-6880-8003-b6a0-bd4d39a6ea04}, 2025.
\newblock Accessed: 2025-06-23.

\end{thebibliography}

Correspondence: \href{mailto:agbanwajamal03@gmail.com}{agbanwajamal03@gmail.com}.

\end{document}